# A restricted additive Vanka smoother for geometric multigrid

S. Saberi[*]    G. Meschke[†]    A. Vogel[†]

**Abstract**

The solution of saddle-point problems, such as the Stokes equations, is a challenging task, especially in large-scale problems. Multigrid methods are one of the most efficient solvers for such systems of equations and can achieve convergence rates independent of the problem size. The smoother is a crucial component of multigrid methods and significantly affects its overall efficiency. We propose a Vanka-type smoother that we refer to as Restricted Additive Vanka and investigate its convergence in the context of adaptive geometric multigrid methods for the Stokes equations. The proposed smoother has the advantage of being an additive method and provides favorable properties in terms of algorithmic complexity, scalability and applicability to high-performance computing. We compare the performance of the smoother with two variants of the classical Vanka smoother using numerical benchmarks for the Stokes problem. We find that the restricted additive smoother achieves comparable convergence rates to the classical multiplicative Vanka smoother while being computationally less expensive per iteration, which results in faster solution runtimes.

## 1 Introduction

Saddle-point systems arise in many applications in science and engineering, such as computational fluid dynamics, structural mechanics, economics, etc. The efficient solution of such systems is a challenging task. Uzawa methods [33, 18, 9], Krylov subspace solvers [40] and multigrid methods [52, 23, 53] are the most popular iterative solvers for saddle-point systems. A comparison of the performance of some of these methods is given in [17, 38, 30]. Methods such as GMRES and Bi-CGSTAB [40], the conjugate residual method in [39, 48] and the conjugate gradient method in [7] fall in the category of Krylov subspace solvers. Appropriate preconditioning is necessary to achieve scalable convergence in such methods. In the case of Stokes and Navier-Stokes systems, it is shown that the preconditioning operator must be spectrally equivalent to the discrete Laplacian terms [48]. Block factorization or Schur complement methods, including pressure correction schemes such as SIMPLE [36] and approximate commutator schemes, such as the pressure convection-diffusion method [47, 26] and the least square commutator method [15] are one class of preconditioners for saddle-point problems, see also [16]. Another class of preconditioners are based on domain decomposition methods such as Schwarz methods [11, 27, 37, 28, 43], see also [12]. Multigrid methods are also often used as preconditioners in Krylov subspace solvers. In the general case, block factorization schemes try to find an appropriate approximation to the velocity block and the Schur complement block separately, where the coupling between velocity and pressure

---

[*]High Performance Computing in the Engineering Sciences, Ruhr University Bochum, Universitätsstr. 150, 44801 Bochum, Germany

[†]Institute for Structural Mechanics, Ruhr University Bochum, Universitätsstr. 150, 44801 Bochum, Germany

is at least partially ignored. On the other hand, monolithic multigrid approaches consider the entire system, preserving the coupling between velocity and pressure.

Multigrid methods have been shown to be among the most efficient solvers for saddle-point problems [17, 44], resulting in convergence rates independent of the mesh size. The smoother is perhaps the most crucial component of multigrid methods. A variety of smoothers have been used for saddle-point systems. A class of smoothers known as Vanka-type smoothers based on the solution of local saddle-point problems was proposed in [51]. A smoother based on incomplete factorization (ILU) was used in [53]. Smoothers based on some pressure correction methods were studied in [21]. The Braess-Sarazin smoother proposed in [6] is a pressure correction scheme and can be derived from the SIMPLE method; however, it differs from SIMPLE-type smoothers, for example, in the fact that the error in Braess-Sarazin smoother depends only on the velocity, while the error in SIMPLE-type smoothers depends only on the pressure. Unlike Vanka-type smoothers, the Braess-Sarazin smoother requires the solution of a global saddle-point problem and updates all degrees of freedom at once. A comparison between Vanka-type smoothers and Braess-Sarazin-type smoothers is contained in [25]. Multiplicative and restricted additive Schwarz smoothers were used for immersed boundary methods in [22] and [5], respectively.

The original Vanka smoother was developed as an overlapping block Gauss-Seidel scheme for the finite difference solution of Navier-Stokes equations [51], which we refer to as the Vanka smoother or the classical or standard Vanka smoother in this work. Local Fourier analysis of the Vanka smoother was studied in [49, 32]. A convergence analysis of the Vanka smoother for the mixed finite element formulation was presented in [35] for the mixed Poisson problem and in [30] for the Stokes equations. The convergence of some Vanka-type smoothers for the Stokes and Navier-Stokes equations was proved in [34]. A diagonal variant of the Vanka smoother was used in [50, 4]. While Vanka-type smoothers are well studied, their application has been typically limited to the multiplicative case, see e.g., [8, 24, 1, 14]. The additive Vanka smoother was studied in [45] and shown to be equivalent to an inexact Uzawa iteration under certain conditions. Recently, local Fourier analysis was used to study a two-level multigrid method with the additive Vanka smoother in [19].

Although the classical Vanka smoother has been shown to be effective for a variety of problems, its multiplicative nature entails some inhibiting drawbacks. First, the multiplicative application of the smoother requires an update to the residual for every subdomain, which is a computationally expensive operation. Second, parallelization of the multiplicative smoother is a non-trivial task on both shared-memory and distributed-memory systems, requires substantial communication and the smoother operator is typically altered to an inexact representation that does not match the original operator in serial, which may compromise the scalability of the smoother from the perspective of high-performance computing for high-fidelity solutions.

Interpreting Vanka-type smoothers from the perspective of Schwarz domain decomposition methods [46, 31, 20] and considering the saddle-point structure of the Stokes equations, we propose a Vanka-type smoother in this work, which can be interpreted as a restricted additive Schwarz method [13]. Our main motivation is to develop an effective smoother for such systems that allows additive application. The main idea of the smoother is to use pressure-oriented subdomains with a modification to the subdomain prolongation operator. We focus on monolithic adaptive geometric multigrid methods and use a Q1-Q1 discretization of the Stokes equations as model problem. The main contributions of this work can be summarized as follows:

- We propose a smoother for saddle-point problems in the context of the Stokes equations with favorable properties in terms of complexity and scalability due to its additive nature

- We investigate the performance of the proposed smoother in the context of multilevel monolithic adaptive geometric multigrid for the Stokes problem

- We compare the performance of the new smoother to other variants of the Vanka smoother using numerical benchmarks

The remainder of this paper is organized as follows. An adaptive geometric multigrid method is formulated for the Stokes equations, adopted as the model problem, in section 2. The proposed smoother is developed in section 3. We analyze the performance of the proposed smoother in comparison with the classical Vanka smoother and the additive Vanka smoother by means of selected benchmark examples in section 4. Finally, we draw some conclusions in section 5.

## 2 Geometric multigrid for the model problem

We start by formulating the finite element discretization of the stationary Stokes equations, which we will use as a model problem. The strong form of the stationary Stokes equations is given by

$$\begin{aligned} -\eta \nabla^2 \boldsymbol{u} + \nabla p &= \boldsymbol{f} \quad \text{in } \Omega, \\ \nabla \cdot \boldsymbol{u} &= 0 \quad \text{in } \Omega, \end{aligned} \quad (1)$$

where $\boldsymbol{u}$ is a vector-valued function representing the fluid velocity, $p$ is a scalar function representing the fluid pressure, $\boldsymbol{f}$ denotes a body force and $\eta$ is the fluid viscosity. Equation eq. (1) in conjunction with the following boundary conditions define the Stokes boundary value problem.

$$\begin{aligned} \boldsymbol{u} &= \boldsymbol{w} \quad \text{on } \Gamma_D \subset \partial\Omega, \\ \eta \frac{\partial \boldsymbol{u}}{\partial \boldsymbol{n}} - \boldsymbol{n} p &= \boldsymbol{h} \quad \text{on } \Gamma_N = \partial\Omega \setminus \Gamma_D, \end{aligned} \quad (2)$$

where $\Omega$ denotes the computational domain, $\partial\Omega$ denotes the boundary of the domain, $\Gamma_D$ and $\Gamma_N$ denote the Dirichlet and Neumann parts of the boundary, respectively and $\boldsymbol{n}$ is the outer normal vector to the boundary.

The infinite-dimensional weak form of the Stokes equations is obtained by multiplying the strong form by appropriate test functions, integrating over the domain and using integration by parts to transfer the derivatives of the pressure trial function to the velocity test functions to reduce the continuity requirements:

Find $(\boldsymbol{u}, p) \in (\boldsymbol{V}, Q)$

$$\begin{aligned} \int_\Omega \eta \nabla \boldsymbol{v} : \nabla \boldsymbol{u} \, d\boldsymbol{x} - \int_{\Gamma_N} \boldsymbol{v} \cdot (\boldsymbol{n} \cdot \eta \nabla \boldsymbol{u}) \, d\boldsymbol{s} & \\ - \int_\Omega \nabla \cdot \boldsymbol{v} p \, d\boldsymbol{x} + \int_{\Gamma_N} \boldsymbol{v} \cdot (\boldsymbol{n} p) \, d\boldsymbol{s} &= \int_\Omega \boldsymbol{v} \boldsymbol{f} \, d\boldsymbol{x}, \quad \forall \boldsymbol{v} \in \boldsymbol{V}_0, \\ \int_\Omega q \nabla \cdot \boldsymbol{u} \, d\boldsymbol{x} &= 0, \quad \forall q \in Q, \end{aligned} \quad (3)$$

where $\boldsymbol{v}$ are the vector-valued velocity test functions, $q$ is the scalar pressure test function and

$$\begin{aligned} \boldsymbol{V} &:= \{\boldsymbol{u} \in H^1(\Omega)^d \mid \boldsymbol{u} = \boldsymbol{w} \text{ on } \Gamma_D\}, \\ \boldsymbol{V}_0 &:= \{\boldsymbol{v} \in H^1(\Omega)^d \mid \boldsymbol{v} = \boldsymbol{0} \text{ on } \Gamma_D\}, \\ Q &:= \{q \in L^2(\Omega)^d\}, \end{aligned} \quad (4)$$

where $H$ and $L$ denote the Sobolev-Hilbert and Lebesgue vector spaces, respectively, and $d$ denotes the dimension.

Let $\Omega_h$ be an appropriate discretization of $\Omega$. Combining the boundary terms in eq. (3) and using a shorthand notation, the discretized weak form can be written as

$$\begin{aligned} (\eta \nabla \boldsymbol{v}_h, \nabla \boldsymbol{u}_h)_{\Omega_h} - (\nabla \cdot \boldsymbol{v}_h, p_h)_{\Omega_h} - (q_h, \nabla \cdot \boldsymbol{u}_h)_{\Omega_h} &= \\ (\boldsymbol{v}_h, \boldsymbol{f}_h)_{\Omega_h} + (\boldsymbol{v}_h, \boldsymbol{n} \cdot (\eta \nabla \boldsymbol{u}_h - p_h \boldsymbol{I}))_{\Gamma_{N_h}}, \end{aligned} \quad (5)$$

where $\boldsymbol{n} \cdot (\eta \nabla \boldsymbol{u}_h - p_h \boldsymbol{I}) = \boldsymbol{h}_h$ over $\Gamma_{N_h}$, see Equation eq. (2). The discretized weak formulation leads to a system of equations of the form

$$\begin{bmatrix} \boldsymbol{A} & \boldsymbol{B} \\ \boldsymbol{B}^T & \boldsymbol{C} \end{bmatrix} \begin{bmatrix} \boldsymbol{u} \\ \boldsymbol{p} \end{bmatrix} = \begin{bmatrix} \boldsymbol{f} \\ \boldsymbol{0} \end{bmatrix}, \tag{6}$$

where $\boldsymbol{A}$ and $\boldsymbol{B}$ are defined by the following bilinear forms

$$\begin{aligned} \boldsymbol{A} &:= (\eta \nabla \boldsymbol{v}_h, \nabla \boldsymbol{u}_h)_{\Omega_h}, \\ \boldsymbol{B} &:= -(\nabla \cdot \boldsymbol{v}_h, p_h)_{\Omega_h}, \end{aligned} \tag{7}$$

and $\boldsymbol{C}$ is zero for a set of conforming spaces that satisfy the inf-sup condition or the stabilization matrix for non-conforming discretizations. We consider a stabilized $Q1$-$Q1$ discretization in this work, which is non-conforming for the Stokes equations, see, e.g., [29]. The following stabilization matrix is used

$$\boldsymbol{C} := -\beta \sum_{\Omega_i^e \in \Omega_h} h_{\Omega_i^e}^2 (\nabla \boldsymbol{v}_h, \nabla \boldsymbol{u}_h)_{\Omega_i^e}, \tag{8}$$

where $\beta$ is a sufficiently large constant, $\Omega_i^e$ is the domain of an element in the discretized domain $\Omega_h$ and $h_{\Omega_i^e}$ is the diameter of $\Omega_i^e$.

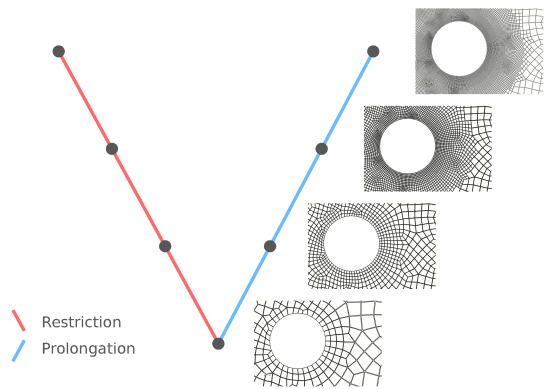

Figure 1: The multigrid $V$-cycle with four levels of grid hierarchy and a combination of uniform and adaptive mesh refinement around a cylinder

We now discuss the construction of a geometric multigrid method, which we will employ for the solution of the system of linear equations given in eq. (6). We consider monolithic multilevel multigrid methods for the solution of the model problem in the primitive variables. The idea of a two-level geometric multigrid method with a fine grid $\Omega_h^l$ and a coarse grid $\Omega_h^{l-1}$ is to smooth highly oscillatory errors on the fine grid and correct slowly oscillatory errors from the coarse grid. This idea can be extended to a hierarchy of grids to form multilevel multigrid methods, see, e.g., [23]. We consider adaptive geometric multigrid with nested spaces. Given a fine grid $\Omega_h^l$ and $l-1$ levels of coarse grids, $\{\Omega_h^1 \ldots \Omega_h^l\}$, a multigrid $V$-cycle can be illustrated as in section 2.

We employ space tree data structures [10] to handle the spatial discretization, where an unstructured mesh consisting of quadrilateral (2D) or hexahedral (3D) elements is mapped onto a set of trees. This mesh, commonly referred to as the macro-mesh, forms the coarsest possible grid. It is then possible to perform an arbitrary number of uniform and adaptive refinement steps. A series of grids with a combination of uniform and adaptive mesh refinement around a cylinder is shown in section 2. The grid hierarchy for the multigrid method and the transfer operators between adjacent grid levels are generated from the refined space tree structure, see, e.g., [42, 41]. A 2:1 balance between neighboring child nodes is maintained, and hanging nodes, which are in general present in adaptively refined grids, are handled as constraints and removed from the global system of equations. In addition to the grid hierarchy and transfer operators, a coarse-grid solver and appropriate smoothers are required for the construction of the multigrid method. We use LU factorization to solve the coarse system to machine accuracy and develop the smoother operator in section 3.

## 3 Restricted additive Vanka smoother

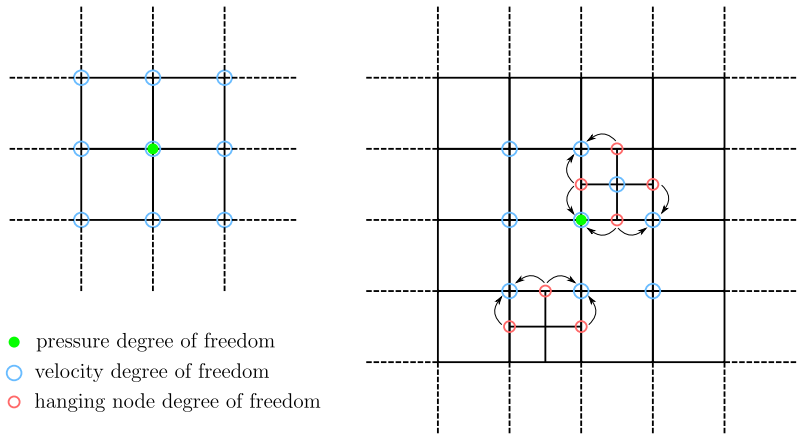

Figure 2: The classical Vanka subdomains for a Q1-Q1 discretization without hanging nodes (left) and with hanging nodes (right). The arrows show how hanging node degrees of freedom can be implicitly included in the subdomain

We discuss Vanka smoothers from the perspective of Schwarz domain decomposition methods. At the discretized level, consider a system of equations of the form

$$\boldsymbol{L}\boldsymbol{x} = \boldsymbol{b}, \tag{9}$$

where $\boldsymbol{L} := \begin{pmatrix} \boldsymbol{A} & \boldsymbol{B} \\ \boldsymbol{B}^T & \boldsymbol{C} \end{pmatrix}$, $\boldsymbol{x} := \begin{pmatrix} \boldsymbol{u} \\ \boldsymbol{p} \end{pmatrix}$ and $\boldsymbol{b} := \begin{pmatrix} \boldsymbol{f} \\ \boldsymbol{0} \end{pmatrix}$ are defined according to Equation eq. (6), and the smoother operator $\boldsymbol{S}$ is defined as

$$\boldsymbol{x}^{k+1} = \boldsymbol{x}^k + \boldsymbol{S}(\boldsymbol{b} - \boldsymbol{L}\boldsymbol{x}^k), \tag{10}$$

where $k$ denotes the iteration step.

Let $\mathcal{S}_x = \{x_1, \ldots, x_n\}$, $\mathcal{S}_p = \{p_1, \ldots, p_{n_p}\}$ and $\mathcal{S}_u = \{u_1, \ldots, u_{n_u}\}$ be the sets of all degrees of freedom and the pressure and velocity degrees of freedom, respectively, where $n$ denotes the total number of degrees of freedom and $n_p$ and $n_u$ denote the number of pressure and velocity degrees of freedom, respectively. Let $V$ be the space spanned by all degrees of freedom over $\mathbb{R}^n$. The classical Vanka smoother for the finite element discretization of Stokes and Navier-Stokes problems is then typically defined as a set of $n_p$ subdomains $\mathcal{S}_{\Omega_{\mathrm{sd}}} = \{\Omega_{\mathrm{sd},1}, \ldots, \Omega_{\mathrm{sd},n_p}\}$, where subdomain $\Omega_{\mathrm{sd},i}$ spans a subspace $V_i \subset V$ over $\mathbb{R}^{n_i}$ that consists of pressure DoF $p_i$ and all the velocity DoFs that are connected to it. By connection, the coupling between velocity and pressure is meant, which topologically corresponds to the vertices of the elements that include the vertex corresponding to $p_i$. Algebraically, it is equivalent to the degrees of freedom corresponding to the nonzero entries in the $i$th row of $B^T$ as in eq. (6) plus the $i$-th pressure degree of freedom. A sample subdomain for a Q1-Q1 discretization is illustrated in section 3. Note that the hanging node degrees of freedom that are implicitly included in the subdomains do not impose an additional computational cost since hanging node constraints are generally resolved prior to the application of the smoother.

Starting with the multiplicative Vanka smoother, the subdomain corrections are then applied in a multiplicative fashion. This decomposition can be understood as an overlapping multiplicative Schwarz method. The smoother operator can be written as

$$\boldsymbol{S}_{\text{MV}} = \prod_{i=1}^{n_{\text{sd}}}(\boldsymbol{R}_{\text{sd},i}^T \boldsymbol{\omega}_i \boldsymbol{L}_i^{-1} \boldsymbol{R}_{\text{sd},i}), \quad (11)$$

where $\boldsymbol{R}_{\text{sd},i} : \mathbb{R}^n \to \mathbb{R}^{n_i}$ is the subdomain restriction operator, $\boldsymbol{L}_i : \mathbb{R}^{n_i} \to \mathbb{R}^{n_i}$ given by $\boldsymbol{L}_i := \boldsymbol{R}_{\text{sd},i} \boldsymbol{L} \boldsymbol{R}_{\text{sd},i}^T$ is the local system of $\Omega_{\text{sd},i}$ and $\boldsymbol{\omega}_i : \mathbb{R}^{n_i} \to \mathbb{R}^{n_i}$ is, in the general case, a diagonal weighting matrix. In general, $n_{\text{sd}}$ is the number of subdomains, which equals to $n_p$ for the classical Vanka method. The subdomain restriction operator $\boldsymbol{R}_{\text{sd},i}$ essentially extracts the set of degrees of freedom $\mathcal{S}_{x,\Omega_{\text{sd},i}} := \{x_j \in \mathcal{S}_x \mid x_j \in \Omega_{\text{sd},i}\}$, and its transpose $\boldsymbol{R}_{\text{sd},i}^T : \mathbb{R}^{n_i} \to \mathbb{R}^n$ acts as an injection operator that prolongates a vector from the subspace $V_i$ of the subdomain to the global space $V$ by padding it with zeros.

The additive version of the smoother operator can be written as

$$\boldsymbol{S}_{\text{AV}} = \sum_{i=1}^{n_{\text{sd}}}(\boldsymbol{R}_{\text{sd},i}^T \boldsymbol{\omega}_i \boldsymbol{L}_i^{-1} \boldsymbol{R}_{\text{sd},i}), \quad (12)$$

where subdomain corrections are applied simultaneously for all subdomains.

We develop an additive smoother by introducing a modification to the prolongation operator in eq. (12). $\boldsymbol{R}_{\text{sd},i}$ can be considered as consisting of those rows of the identity matrix $\boldsymbol{I}^n \in \mathbb{R}^{n \times n}$ that correspond to the degrees of freedom in $\Omega_{\text{sd},i}$. We define $\tilde{\boldsymbol{R}}_{\text{sd},i} : \mathbb{R}^n \to \mathbb{R}^{n_i}$ by letting all rows in $\boldsymbol{R}_{\text{sd},i}$ be zero except for those corresponding to the pressure DoF and the velocity DoFs on the pressure node. In other words, let $\boldsymbol{G}_i := \tilde{\boldsymbol{R}}_{\text{sd},i}^T \tilde{\boldsymbol{R}}_{\text{sd},i} \in \mathbb{R}^{n \times n}$ be a diagonal matrix. Then, the set $\mathcal{S}_{x,\tilde{\Omega}_{\text{sd},i}} := \{x_j \in \mathcal{S}_x \mid (G_i)_{j,j} = 1\}$ only contains the pressure DoF and the velocity DoFs on the pressure node of $\Omega_{\text{sd},i}$. The proposed smoother operator can then be defined as

$$\boldsymbol{S}_{\text{RAV}} = \sum_{i=1}^{n_{\text{sd}}}(\tilde{\boldsymbol{R}}_{\text{sd},i}^T \boldsymbol{\omega}_i \boldsymbol{L}_i^{-1} \boldsymbol{R}_{\text{sd},i}), \quad (13)$$

where $\tilde{\boldsymbol{R}}_{\text{sd},i}^T$ is the modified prolongation operator. This decomposition is analogous to the idea of the restricted additive Schwarz method [13]. Therefore, we refer to it as restricted additive Vanka (RAV). To our best knowledge, it is the first time this smoother is used for the solution of the Stokes problem. We note that the proposed smoother can conceptually be applied to other saddle-point problems.

## 4 Numerical experiments

We consider a number of numerical studies to evaluate the performance of the proposed restricted additive Vanka smoother and compare it to the existing multiplicative and additive Vanka methods. All variants are analyzed as smoothers in the geometric multigrid method. Two benchmarks, namely driven cavity and cylinder flow are considered, in which Stokes flow with appropriate boundary conditions is solved. The fluid viscosity is $\eta = 10^{-3} m^2/s$.

Geometric multigrid method is used as a solver in all examples. The iterations are terminated when the initial residual is reduced by a factor of $10^8$. The Vanka smoother with a predefined number of smoothing steps is set up on all levels of the hierarchy. Always, the same number of pre- and post-smoothing steps is applied. The system of equations is assembled on all grids. A damping factor of 0.1 is used for the additive Vanka smoother and a damping factor of 0.66 is used for the restricted additive Vanka and multiplicative Vanka smoothers, which we found were necessary to achieve robust convergence in all considered cases. We note that the geometric multigrid method can also be used as a preconditioner in a Krylov subspace solver, which typically leads to an improvement in performance. Nevertheless, we consider geometric multigrid as a standalone fixed-point iterator in order to better isolate the behavior of the investigated variants of smoothers.

In each benchmark, we analyze the convergence of the geometric multigrid solver through a mesh study, where a set of $n_g$ grids $\Omega_h^i, i = 2, \ldots, n_g$ is generated via uniform and adaptive refinement from a coarse grid $\Omega_h^1$ to form a hierarchy of nested spaces. $M_{\mathfrak{S}}^i, i = 2, \ldots, n_g$ denotes the multigrid operator of a problem with $\Omega_h^i$ as the fine grid and $i$ levels of grid hierarchy, i.e., all problems use $\Omega_h^1$ as the coarse grid. $\mathfrak{S}$ denotes the smoother type, i.e., one of AV (additive Vanka), RAV (restricted additive Vanka) or MV (multiplicative Vanka). The criteria for mesh refinement are problem dependent. A 2:1 balance is maintained on all grids, which implies that some regions of the mesh that lie outside of the targeted refinement zone may also be refined. Three steps of pre-smoothing and three steps of post-smoothing are used for all problems in the convergence study. In the second benchmark (cylinder flow), we also investigate the effectiveness of each smoother as a function of the smoothing steps. We further study the performance of the geometric multigrid solver by analyzing the computational cost of the different smoothers. All examples are conducted using an in-house C++ code with `p4est` [10] and `PETSc` [3] libraries for mesh handling and parts of the linear algebra, respectively, and run on a single core of an Intel Xeon Skylake Gold 6148 CPU with 192 GB of main memory. `Paraview` [2] is used for postprocessing the results.

## 4.1 Benchmark 1: Driven cavity problem

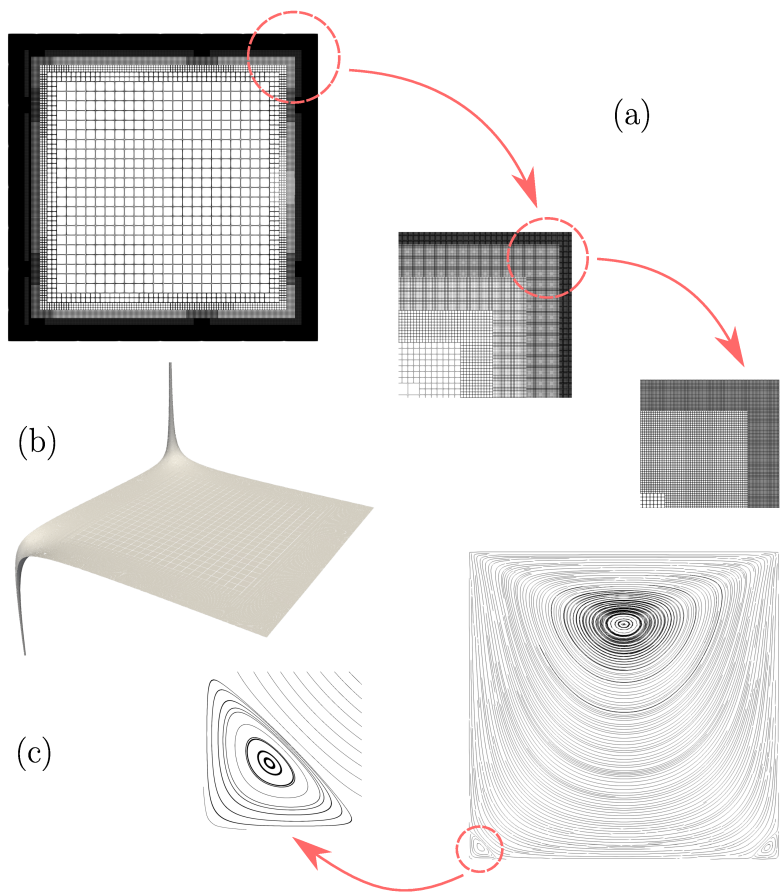

Figure 3: Driven cavity benchmark: (a) Adaptive refinement towards the boundaries, (b) pressure distribution and (c) streamlines

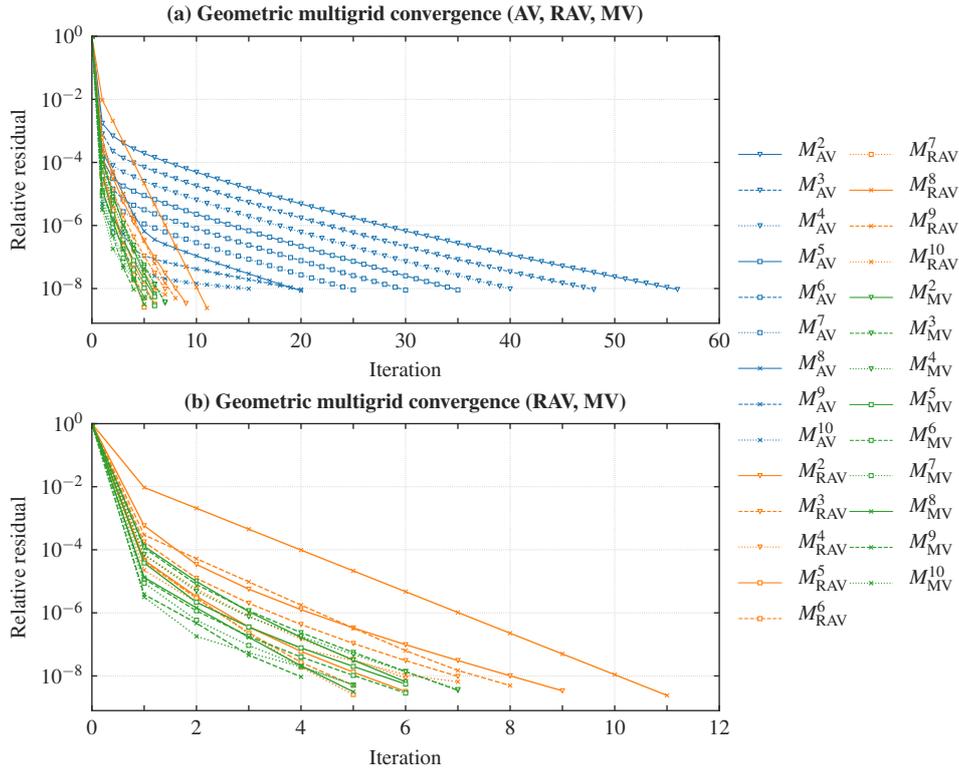

Figure 4: Driven cavity benchmark: Convergence study of the geometric multigrid solver. (a) Multiplicative Vanka, additive Vanka and restricted additive Vanka smoothers and (b) a close-up look at the multiplicative Vanka and the proposed restricted additive Vanka smoothers

The driven cavity problem is an enclosed flow in a square domain, where no-slip wall boundary conditions are imposed on the left, right and bottom sides, and a zero vertical velocity and a nonzero horizontal velocity are imposed on the top side. We use the unit square as the computation domain and impose a horizontal velocity of $u_x = 1.0$ $m/s$ on the top boundary. Starting from the unit square, we use a grid with three levels of uniform refinement as the coarse grid and generate a total of ten grids by adaptively refining the domain towards the boundaries. At each refinement step, the refinement area becomes smaller. The refinement strategy is illustrated in section 4.1(a). The grid hierarchy used for the convergence study is summarized in table 1. The convergence of the geometric multigrid solver is shown in section 4.1 for different smoothers. It can be seen that the additive Vanka smoother (AV) consistently requires the most number of iterations for convergence, while the number of required iterations for the restricted additive Vanka (RAV) and the multiplicative Vanka (MV) is substantially smaller.

The average reduction factor, runtime per iteration and total runtime of the solver for each investigated smoother is shown in section 4.1. While the restricted additive smoother and the multiplicative smoother achieve comparable reduction factors, a salient gap between the additive smoother and the multiplicative and restricted additive smoothers is maintained. On the other hand, the restricted additive smoother and the additive smoother are considerably cheaper per iteration in comparison with the multiplicative smoother, especially on finer grids. The combination of these metrics determines the total runtime of the solver. The total runtime of multiplicative and additive smoothers are similar, with the multiplicative smoother being marginally faster on most grids. It shows how the multiplicative smoother enjoys higher effectiveness but suffers from its higher cost per iteration. The opposite can be said about the additive smoother. On the other hand, the proposed restricted additive smoother is the cheapest of the three methods per iteration and achieves reduction rates that are similar to the multiplicative smoother. As a consequence, it consistently achieves the fastest total runtime. section 4.1(c) demonstrates that the runtime of the proposed restricted additive smoother is up to more than ten times smaller compared to the other variants in this benchmark.

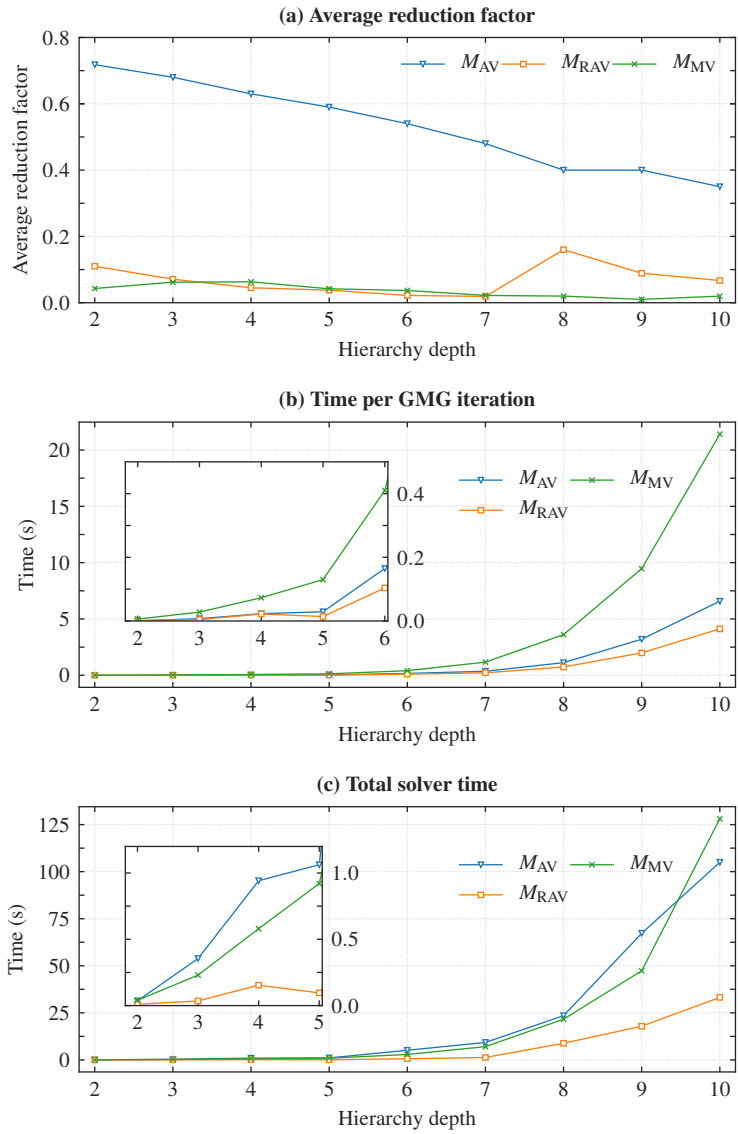

Figure 5: Driven cavity benchmark: Performance of the geometric multigrid solver w.r.t. the grid hierarchy: (a) Average reduction factor of the solver, (b) runtime per iteration and (c) the total runtime of the solver

Table 1: Driven cavity benchmark: The grid hierarchy. $n_e$ and $n_{DoF}$ denote the number of elements and degrees of freedom, respectively

| Grid | $n_e$ | $n_{\text{DoF}}$ |
|---|---|---|
| $\Omega_h^{10}$ | 1,020,508 | 3,063,075 |
| $\Omega_h^{9}$ | 606,184 | 1,819,899 |
| $\Omega_h^{8}$ | 236,788 | 711,339 |
| $\Omega_h^{7}$ | 79,648 | 239,595 |
| $\Omega_h^{6}$ | 25,612 | 77,259 |
| $\Omega_h^{5}$ | 7,672 | 23,283 |
| $\Omega_h^{4}$ | 2,296 | 7,059 |
| $\Omega_h^{3}$ | 676 | 2,139 |
| $\Omega_h^{2}$ | 208 | 699 |
| $\Omega_h^{1}$ | 64 | 243 |

Table 2: Cylinder flow benchmark: The grid hierarchy. $n_e$ and $n_{DoF}$ denote the number of elements and degrees of freedom, respectively

| Grid | $n_e$ | $n_{\text{DoF}}$ |
|---|---|---|
| $\Omega_h^{7}$ | 556,508 | 1,670,814 |
| $\Omega_h^{6}$ | 300,833 | 903,648 |
| $\Omega_h^{5}$ | 96,800 | 291,564 |
| $\Omega_h^{4}$ | 31,082 | 94,464 |
| $\Omega_h^{3}$ | 12,128 | 37,392 |
| $\Omega_h^{2}$ | 3,032 | 9,600 |
| $\Omega_h^{1}$ | 758 | 2,526 |

## 4.2 Benchmark 2: Channel flow with a cylindrical obstacle

The cylinder flow benchmark is an inflow / outflow problem inside a channel characterized by an embedded cylindrical obstacle located in the vicinity of the inflow boundary, see, e.g., [44]. The channel dimensions are $2.2\ m \times 0.41\ m$. The cylindrical obstacle is located at $(0.2\ m, 0.2\ m)$ and has a radius of $0.05\ m$. The inflow boundary condition is given by

$$
\begin{aligned}
u_x(0,y) &= \frac{4\bar{u}y(H-y)}{H^2} \quad y \in [0, H], \\
u_y(0,y) &= 0 \quad y \in [0, H],
\end{aligned}
\tag{14}
$$

where $\bar{u} = 0.3\ m/s$ is the mean horizontal velocity, and $H$ is the channel hight. A homogeneous Neumann condition is imposed on the outflow, and no-slip wall conditions are imposed on the top and bottom boundaries of the channel as well as along the cylindrical body embedded in the flow.

A 2D finite element mesh with 758 Q1-Q1 elements is used as the course grid, which is adaptively refined towards the boundaries to construct the grid hierarchy. The refinement region becomes progressively smaller at each refinement step as illustrated in section 4.2(a).

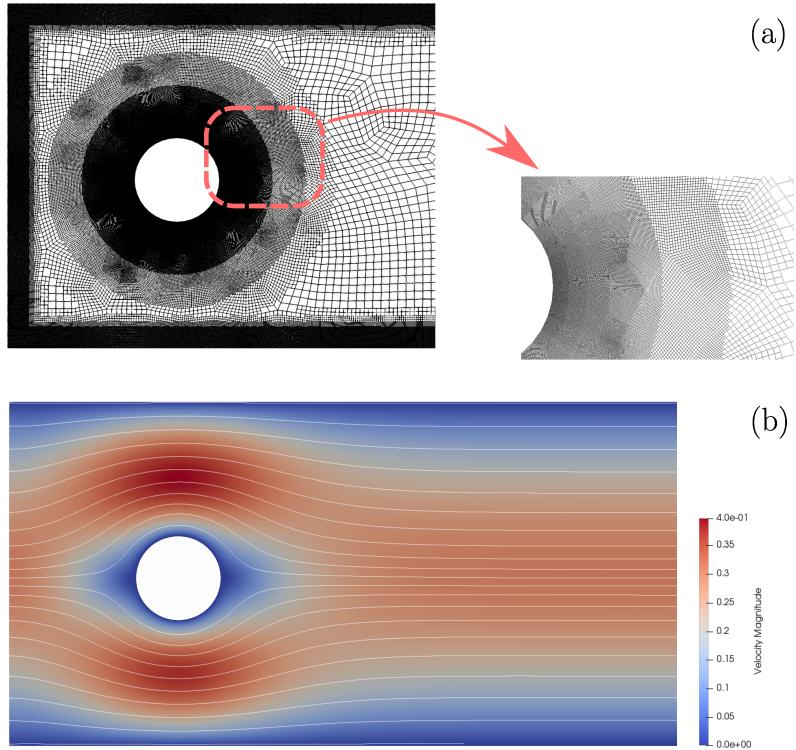

Figure 6: Cylinder flow benchmark: (a) Adaptive refinement of the computational mesh towards the boundaries and (b) contours of the velocity magnitude and streamlines. Note that only the left portion of the channel is shown

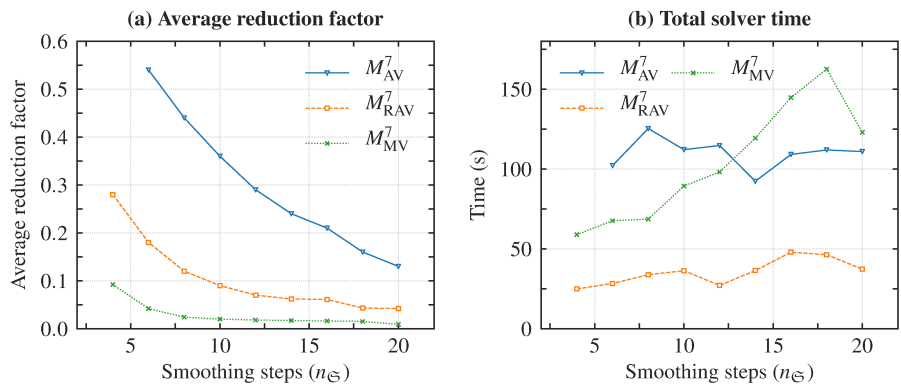

Figure 7: Cylinder flow benchmark: Influence of the number of smoothing steps on the convergence of the geometric multigrid solver. (a) Average reduction factor and (b) total solver runtime. The number of smoothing steps is the sum of $n_{\mathfrak{S}}$ pre-smoothing and $n_{\mathfrak{S}}$ post-smoothing steps. Note that the additive Vanka smoother requires a minimum of three pre- and post-smoothing steps to achieve convergence

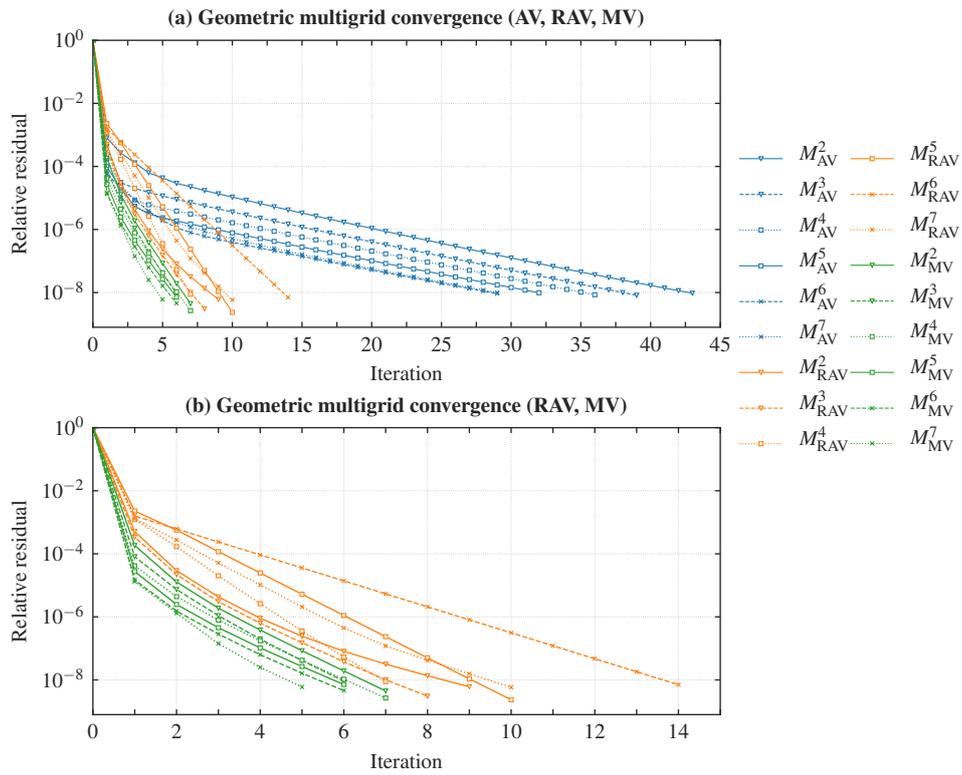

Figure 8: Cylinder flow benchmark: Convergence study of the geometric multigrid solver. (a) Multiplicative Vanka, additive Vanka and restricted additive Vanka smoothers and (b) a close-up look at the multiplicative Vanka and the proposed restricted additive Vanka smoothers

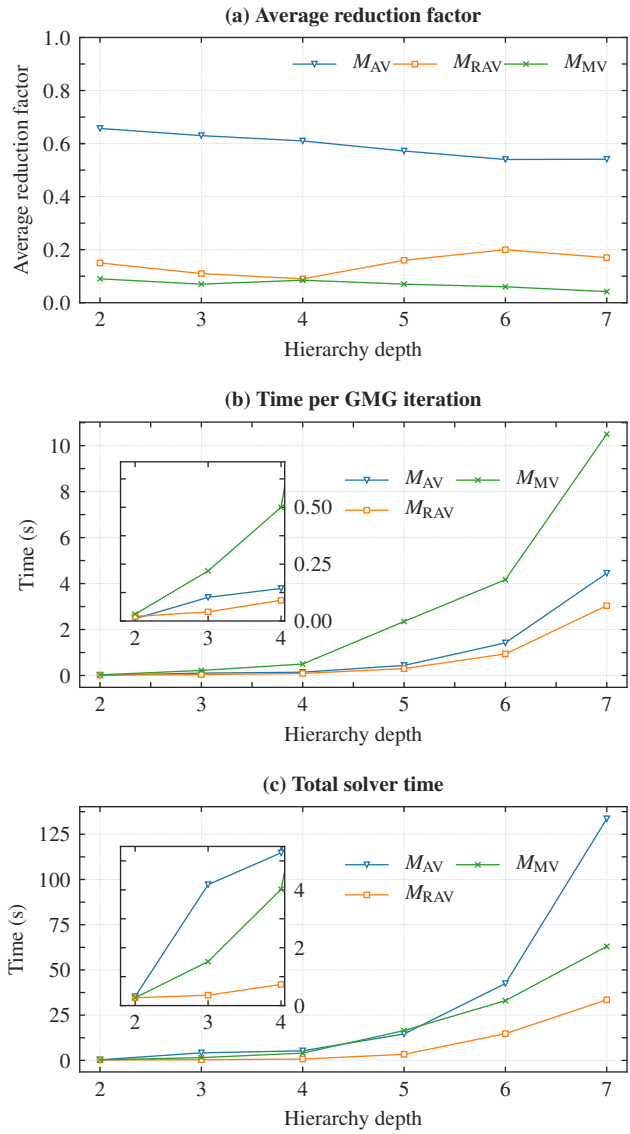

Figure 9: Cylinder flow benchmark: Performance of the geometric multigrid solver w.r.t. the grid hierarchy. (a) Average reduction factor of the solver, (b) runtime per iteration and (c) total runtime of the solver

The number of elements and degrees of freedom for the resulting grid hierarchy is summarized in table 2. The contours plot of the velocity magnitude and the streamlines are shown in section 4.2(b).

We first analyze the effectiveness of the geometric multigrid solver as a function of the number of smoothing steps. Therefore, the finest problem with the deepest grid hierarchy, i.e., $M^7$, is solved using an increasingly larger number of smoothing steps. The average reduction factor as well as the total runtime of the solver are shown in section 4.2. For a given number of smoothing steps, the multiplicative Vanka smoother is consistently the most effective, followed by the restricted additive Vanka smoother, and the additive Vanka smoother is the least effective of the three variants. However, the gap between the multiplicative and restricted additive Vanka smoothers is small and becomes narrower at higher smoothing steps. It can also be seen that, while the effectiveness of the additive smoother drastically changes with the number of smoothing steps, the multiplicative and restricted additive smoothers reach a plateau after a relatively small number of smoothing steps. The latter behavior is in general desirable and indicates that the geometric multigrid method is robust.

In practice, the total solver runtime is typically a more important metric than the shear reduction capability of the solver per iteration, which provides a different analytical landscape, in whose light it can now be seen that the restricted additive smoother needs by far the lowest computation time, followed by the multiplicative and additive smoothers, respectively. The runtime is reduced approximately by a factor of 3 compared to the additive smoother and a factor ranging between 2 and 5 compared to the multiplicative smoother. While the runtime of the additive and restricted additive smoothers remains roughly constant with the number of smoothing steps, the runtime of the multiplicative smoother consistently increases. This indicates, that the additive smoothers are more cost effective. The computational cost of the smoothers is discussed in more detail in section 4.3.

We study the convergence of the geometric multigrid solver on the grid hierarchy next. As shown in section 4.2, the additive smoother requires considerably more iterations as compared to the multiplicative and restricted additive smoothers on all grid levels. On the other hand, the multiplicative and restricted additive smoothers require a similar number of iterations. The performance of the solver is shown in section 4.2. It can be seen that the slight increase in the number of required iterations in the restricted additive smoother compared to the multiplicative smoother is offset by its cheaper cost per iteration. Therefore, the proposed restricted additive smoother overall leads to the fastest solver runtime. Analogous to the restricted additive smoother, the additive smoother has a low cost per iteration; however, its total runtime is typically larger than the multiplicative smoother because of its iteration count.

## 4.3 Computational cost

We discuss the computational cost of the smoothers in this subsection. The only difference between the smoother operators $S_{\text{MV}}$ and $S_{\text{AV}}$ is the fact that the former requires that the residual is updated after the application of each subdomain $\Omega_{\text{sd},i}$, which in general requires the computation of

$$x_m = b_m - L_{m,k} x_k, \ k = 1, \ldots, n \quad \forall \, \{m \mid x_m \in \Omega_{\text{sd},i}\}, \tag{15}$$

whereas the latter updates all degrees of freedom only once. The difference is more significant in larger problems with deeper grid hierarchies as shown in section 4.1(b) and section 4.2(b).

The restricted additive smoother according to eq. (13) differs from the additive smoother according to eq. (12) only in the prolongation operator. The inverse of the local subdomain system $\boldsymbol{L}_i^{-1} \in \mathbb{R}^{n_i \times n_i}$ is in general a dense matrix; therefore, it has to be computed also for the restricted additive smoother. However, the application of the subdomain inverse can be limited to the space spanned by $\tilde{\Omega}_{\text{sd},i}$ over $\mathbb{R}^{\tilde{n}_i \times \tilde{n}_i}$ with $\tilde{n}_i < n_i$. Therefore, the restricted additive smoother requires fewer floating-point operations when this opportunity is exploited. The results in section 4.1(b) and section 4.2 corroborate this observation. Furthermore, the restricted additive smoother has more favorable memory requirements in case the subdomain inverses are precomputed as it only requires the storage of certain rows.

As shown in section 4.1(a) and section 4.2(a), the higher computational cost of the multiplicative smoother is compensated by its effectiveness, and the total runtime of the multiplicative smoother is often less than the additive smoother despite its larger runtime per iteration. It was also shown that, while the additive smoother promotes the use of more smoothing steps, the multiplicative smoother rewards fewer smoothing steps as the improvement from more steps is overshadowed by the accumulated cost of the smoother. The restricted additive smoother strikes a good balance between computational cost and effectiveness. It is the least expensive of the three methods and yet achieves reduction rates that are comparable to the multiplicative smoother. Consequently, the restricted additive smoother consistently achieves the fastest total solver runtime.

# 5 Conclusions

We proposed a Vanka-type smoother, denoted as restricted additive Vanka (RAV), and investigated its convergence in the context of multilevel monolithic adaptive geometric multigrid methods for solving the Stokes equations. We compared the performance of restricted additive Vanka with multiplicative Vanka and additive Vanka smoothers in the context of the geometric multigrid solver using two numerical benchmarks. The results indicate that the multiplicative Vanka smoother and the additive Vanka smoother are the most and least effective per iteration of geometric multigrid, respectively. On the other hand, the restricted additive Vanka smoother is computationally more favorable and achieves reduction rates that are comparable to the multiplicative Vanka smoother. As a consequence, the proposed restricted additive Vanka smoother consistently achieves the fastest total solver runtime in the investigated benchmarks. Furthermore, the restricted additive Vanka operator has the advantage of being an additive method, which renders it highly attractive from the perspective of parallel computing. The analysis of the scalability of the proposed smoother is a topic for future work.

## Acknowledgments

Financial support was provided by the German Research Foundation (Deutsche Forschungsgemeinschaft, DFG) in the framework of subproject C4 of the collaborative research center SFB 837 Interaction Modeling in Mechanized Tunneling, grant number 77309832. This support is gratefully acknowledged. We also gratefully acknowledge the computing time on the computing cluster of the SFB 837 and the Department of Civil and Environmental Engineering at Ruhr University Bochum, which has been employed for the presented studies.